\documentclass[10pt]{amsart}
\usepackage{import}
\usepackage{./preamble/Pakete}
\usepackage{./preamble/Befehle}

\begin{document}


\title[K-theory Soergel Bimodules]{K-theory Soergel Bimodules}

\author{Jens Niklas Eberhardt}
\address{Bergische Universität Wuppertal\\ 
Gaußstraße 20\\
D-42119 Wuppertal\\
Germany
}
\email{mail@jenseberhardt.com}  
\begin{abstract} 
    We initiate the study of $K$-theory Soergel bimodules---a $K$-theory analog of classical Soergel bimodules. Classical Soergel bimodules can be seen as a completed and infinitesimal version of their new $K$-theoretic analog.
    
    We show that morphisms of $K$-theory Soergel bimodules can be described geometrically in terms of equivariant $K$-theoretic correspondences between Bott--Samelson varieties. We thereby obtain a natural categorification of $K$-theory Soergel bimodules in terms of equivariant coherent sheaves.

    We introduce a formalism of stratified equivariant $K$-motives on varieties with an affine stratification, which is a $K$-theoretric analog of the equivariant derived category of Bernstein--Lunts.
    We show that Bruhat-stratified torus-equivariant $K$-motives on flag varieties can be described in terms of chain complexes of $K$-theory Soergel bimodules.

    Moreover, we propose conjectures regarding an equivariant/monodromic Koszul duality for flag varieties and the quantum $K$-theoretic Satake.
\end{abstract}
\maketitle
\setcounter{tocdepth}{1} 
\tableofcontents

\section{Introduction}
Let $G\supset B\supset T$ be a connected split reductive group with a Borel subgroup $B$ and maximal torus $T$ such that $G'$ is simply connected. Let $N_G(T)/T=\pW\supset \pS$ be the Weyl group with set of simple reflections $\pS$, $X(T)=\Hom_{grp}(T,\Gm)$ the character lattice and $X=G/B$ the flag variety. Let $\Lambda$ be some ring of coefficients.
\subsection{Soergel bimodules} \emph{Soergel bimodules} are graded bimodules over the $T$-equivariant ring of a point $$S=H^{\bullet}_T(\pt)=H^{\bullet}(BT)=\operatorname{Sym}^\bullet(X(T)_\Lambda)$$ that arise as direct summands of  the $T$-equivariant cohomologies of \emph{Bott--Samelson varieties} 
$$\BSvar(s_1,\dots,s_n)=P_{s_1}\times_B\dots\times_BP_{s_n}/B$$
which admit a simple description as iterated tensor products
$$H^\bullet_T(\BSvar(s_1,\dots,s_n))=S\otimes_{S^{s_1}} \dots \otimes_{S^{s_n}}S,$$
see \cite{soergelCombinatoricsHarishChandraBimodules1992}.

If $\Lambda$ is a field of characterstic $0$, indecomposable Soergel bimodules yield the equivariant \emph{intersection cohomology} of \emph{Schubert varieties} in $X$. This is a consequence of the \emph{decomposition theorem} for perverse sheaves, see \cite{beilinsonFaisceauxPervers1982}, and Soergel's \emph{Erweiterungssatz}, see \cite{soergelKategoriePerverseGarben1990}. 

\subsection{$K$-theory Soergel bimodules} Our definition of \emph{$K$-theory Soergel bimodules} follows the simple idea of replacing equivariant cohomology by equivariant $K$-theory: The ring $S$ is replaced by the representation ring of the torus $T$ $$R=R(T)=K_0^T(\pt)=K_0(\Rep(T))=\Lambda[X(T)].$$
The $T$-equivariant $K$-theory of a Bott--Samelson variety is a bimodule over $R$ and can be computed as 
$$K^T_0(\BSvar(s_1,\dots,s_n))=R\otimes_{R^{s_1}} \dots \otimes_{R^{s_n}}R,$$
see \eqref{eq:ktheorybottsamelson}. Consequently, we define $K$-theory Soergel bimodules as direct sums and direct summands of these bimodules. We note that $K$-theory Soergel bimodules were also considered in \cite[Section 8]{dyerRepresentationTheoriesCoxeter1995}.
\subsection{Atiyah--Segal completion theorem} If $\Lambda$ is a field of characterstic $0$, (cohomological) Soergel bimodules can be interpreted as a \emph{completed} or \emph{infinitesimal} version of $K$-theory Soergel bimodules.

Namely, the Atiyah--Segal completion theorem and the Chern character isomorphism exhibit the equivariant cohomology (=cohomology of the Borel construction) as the completion of the genuine equivariant $K$-theory (= $K$-theory of equivariant vector bundles) at the augmentation ideal $I=\ker(K^T_0(\pt)\to K_0(\pt)).$ For a Bott--Samelson variety $\BSvar$ we obtain 
\[\begin{tikzcd}
	{K^T_0(\BSvar)^\wedge_I} & {K_0(ET\times_T\BSvar)} & {\prod_iH^i(ET\times_T\BSvar)=\prod_iH^\bullet_T(\BSvar).}
	\arrow["\sim", from=1-1, to=1-2]
	\arrow["\sim", from=1-2, to=1-3]
\end{tikzcd}\]
 
Geometrically, cohomological and $K$-theory Soergel bimodules can be interpreted as coherent sheaves on the spaces 
$$\mathfrak{t}^\vee\times_{\mathfrak{t}^\vee/W}\mathfrak{t}^\vee\text{ and }T^\vee\times_{T^\vee/W}T^\vee,$$
respectively, where $T^\vee/\Lambda$ denotes the dual torus and $\mathfrak{t}^\vee$ its Lie algebra. Hence, the former arise from the latter by passing to an infinitesimal neighborhood at the identity $1\in T^\vee.$ 
\subsection{Correspondences and Erweiterungssatz} We will show that morphisms of $K$-theory Soergel bimodules admit a geometric description of in terms of $K$-theoretic correspondences between Bott--Samelson varieties.
\begin{theorem*}[Theorem \ref{thm:erweiterungssatzconvolutionktheory}]
    Let $\mathbf{x}, \mathbf{y}$ be sequences of simple reflections. Then convolution induces an isomorphism
    $$\act: G_0^T(\BSvar(\mathbf{x})\times_X\BSvar(\mathbf{y}))\stackrel{\sim}{\to} \Hom_{R\otimes R}(K_0^T(\BSvar(\mathbf{x})),K_0^T(\BSvar(\mathbf{y}))).$$
\end{theorem*}
Here $G_0^T$ denotes $K_0(\Coh_T(-))$ which is the $K$-theoretic analog of Borel--Moore homology. The result is a $K$-theoretic analog of Soergel's Erweiterungssatz, which implies a similar statement for equivariant Borel--Moore homology and cohomology.
Correspondences can be composed via convolution which allows to lift the result to an equivalence of categories between the Karoubi envelope of a category of $K$-theoretric correspondences and $K$-theory Soergel bimodules.
\subsection{$K$-motives on flag varieties} Equivariant cohomology groups $H_T^\bullet$ can be interpreted as extensions in the equivariant derived category of constructible sheaves $\D_T,$ see \cite{bernsteinEquivariantSheavesFunctors1994}. This yields to another construction of (cohomological) Soergel bimodules in terms of equivariant sheaves on the flag variety $X$.

In a modern formulation, this can be described via an equivalence of stable $\infty$-categories
$$\D^{mix}_{T,(B)}(X)\stackrel{\sim}{\to}\Chb(\SBim^\Z_S)$$
between a category of $T$-equivariant mixed sheaves, that are locally constant along Bruhat cells, and the category of chain complexes of graded Soergel bimodules. 

Mixed sheaves $\D^{mix}(X)$ are a graded refinement of the category of constructible sheaves $\D^b(X)$ that can be constructed via mixed Hodge modules or
mixed $\ell$-adic sheaves, see \cite{beilinsonKoszulDualityPatterns1996,hoRevisitingMixedGeometry2022},
and, most satisfyingly, using
mixed motives $\DM(X)$, see for example \cite{soergelPerverseMotivesGraded2018,soergelEquivariantMotivesGeometric2018,eberhardtMixedMotivesGeometric2019, eberhardtIntegralMotivicSheaves2022}.

We will prove a $K$-theoretic analog of this story which provides a third definition of $K$-theory Soergel bimodules. 

Equivariant $K$-theory groups $K_0^T$  can be interpreted as morphisms in the category of \emph{equivariant $K$-motives} $\DK^T$. We will give a definition of this category based on the equivariant stable homotopty category $\SH^T$ constructed in \cite{hoyoisSixOperationsEquivariant2017}. We will see that $\DK^T$ comes equipped with a six functor formalism and behaves very similary to $\D^{mix}_T.$ In particular, we will discuss \emph{affine-stratified} $K$-motives in detail and discuss their formality using weight structures. We will then show:
\begin{theorem*}[Corollary \ref{cor:kmotivesviacomplexesofsoergelbimodules}]
    Let $\Lambda=\Q$. There is an equivalence of stable $\infty$-categories
    $$\DK^T_{(B)}(X)\stackrel{\sim}{\to}\Chb(\SBim_R)$$
    between the category of Bruhat-stratified $T$-equivariant $K$-motives on the flag variety and the category of chain complexes of $K$-theory Soergel bimodules over $R.$
\end{theorem*}
\subsection{Further Directions} This paper should be seen as a starting point to new possible $K$-theoretic approaches to geometric representation theory. We now discuss some of these further directions.
\subsubsection{Categorification of $K$-theory Soergel bimodules} The interpretation of $K$-theory Soergel bimodules and their morphisms in terms of $K$-theory of (fiber products of) Bott--Samelson varieties immediately yields a categorification
\[\begin{tikzcd}[row sep = 1pt]
	{R\otimes_{R^{s_1}}\cdots\otimes_{R^{s_n}}R = K_0^T(\BSvar(\mathbf{x}))} & {\D^b(\Coh_B(\BSvar(\mathbf{x})))} \\
	{\Hom_{R\otimes R}(K_0^T(\BSvar(\mathbf{x})),K_0^T(\BSvar(\mathbf{y})))} & {\D^b(\Coh_B(\BSvar(\mathbf{x})\times_X^{\mathbb{L}}\BSvar(\mathbf{y})))}
	\arrow["{K_0}"', squiggly, from=1-2, to=1-1]
	\arrow["{K_0}"', squiggly, from=2-2, to=2-1]
\end{tikzcd}\]
in terms of the derived category of equivariant coherent sheaves on these spaces. Composition of morphisms is categorified with a convolution formula similar to Fourier--Mukai transformations.
We are hence in the funny situation of having a categorification of a categorification of the Weyl group! We will explore the implications in a future work.
\subsubsection{Diagrammatic Calculus and Algebraic Properties} 
Cohomological Soergel bimodules admit a diagrammatic description that, roughly speaking, describes the relationship between the units and counits induced by the various Frobenius extensions $S^s\subset S$ for $s\in \pS$, see \cite{eliasDiagrammaticsSoergelCategories2011,eliasSoergelCalculus2016,eliasSingularSoergelBimodules2020}.

Very similarly, there are Frobenius extensions $R^s=R(P_s)\subset R=R(T)$ for $s\in \pS$ which arise from parabolic induction. They fulfill similar relationships and it is very imaginable that there is a diagrammatic calculus for $K$-theory Soergel bimodules. For example, there should be a nice diagrammatic basis for their homomorphisms corresponding to the affine strata of the fiber products $\BSvar(\mathbf{x})\times_X\BSvar(\mathbf{y}).$

In this paper we completely ignore any algebraic questions such as a Krull--Schmidt property, uniqueness of indecomposable $K$-theory Soergel bimodules, etc., which are probably best studied using diagrammatics.

\subsubsection{Equivariant/Monodromic duality} 
Koszul duality for flag varieties, see \cite{soergelKategoriePerverseGarben1990, beilinsonKoszulDualityPatterns1996}, is an equivalence of categories between mixed sheaves on a flag variety $X$ and its Langlands dual $X^\vee.$ Equivalently, Koszul duality provides an equivalence of the derived graded principal block of category $\pazocal O$ of a complex reductive Lie algebra and its Langlands dual.

Remarkably, Koszul duality intertwines the Tate-twist and shift functor $(1)[2]$ with the Tate twist $(1)$. This motivated our construction of a non-mixed/ungraded Koszul duality for flag varieties, see \cite{eberhardtMotivesKoszulDualitya},
$$\DK_{(B)}(X)\stackrel{\sim}{\to}\D_{(B)}(X^{\vee}),$$
relating $K$-motives to constructible sheaves: $K$-motives admit a phenomenon called \emph{Bott periodicity} which implies that $(1)[2]$ is the identity functor, while the Tate twist $(1)$ acts trivially on (non-mixed) constructible sheaves.

In the spirit of \cite{bezrukavnikovKoszulDualityKacMoody2014}, this result should have a equivariant/monodromic lift:
\begin{conjecture}[Ungraded, uncompleted equivariant/monodromic Koszul duality] Let $\Lambda=\Q.$ There should be an equivalence of categories
    $$\DK^T_{(B)}(X)\stackrel{\sim}{\to}\D^{b,fg}_{B^\vee\times B^\vee\operatorname{-mon}}(G^{\vee}(\C)),$$
    between Bruhat-stratified $T$-equivariant $K$-motives on a flag variety and Bruhat-stratified $B^\vee\times B^\vee$-monodromic constructible sheaves on the Langlands dual group whose stalks are finitely generated under the fundamental group of $B^\vee\times B^\vee.$
\end{conjecture}
For each maximal ideal $I\subset R$, this conjecture specializes to a Koszul duality between  $I$-twisted equivariant sheaves and $I$-locally finite monodromic sheaves (see \cite{gouttardPerverseMonodromicSheaves2022, lusztigEndoscopyHeckeCategories2020}).
\subsubsection{Quantum $K$-theoretic Satake}
The approach to $K$-theoretic correspondences via $K$-motives developed here in the context of $K$-theory Soergel bimodules should shed new light on Cautis--Kamnitzer's quantum $K$-theoretic Satake, see \cite{cautisQuantumKtheoreticGeometric2018}, which can be reformulated as the following:
\begin{conjecture} There should be an equivalence of categories
    $$\DK^{G\times \Gm}_\text{r}(\operatorname{Gr})\stackrel{\sim}{\to}\D^{b}_{U_q(\mathfrak{g}^\vee)}(\pazocal{O}_q(G^\vee)),$$
    between reduced $G\times \Gm$-equiariant $K$-motives on the affine Grassmannians and $U_q(\mathfrak{g}^\vee))$-equivariant $\pazocal{O}_q(G^\vee)$-modules.
\end{conjecture}
Here reduced $K$-motives $\DK_\text{r}$ should be constructed from $\DK$ by removing the higher $K$-theory of the base point, as defined in the context of $\DM$ in \cite{eberhardtIntegralMotivicSheaves2022}.
\subsection{Structure of the paper} 
In Section \ref{sec:kmotives}, we introduce the formalism of $G$-equivariant $K$-motives $\DK^G$ for diagonalizable groups $G$ and discuss the relation to equivariant $K$-theory and $G$-theory. 

In Section \ref{sec:stratified}, we consider $\Ss$-stratified $G$-equivariant $K$-motives $\DK_\Ss^G(X)$ for varieties $X$ with an affine stratification $\Ss.$ We construct a weight structure and discuss their formality.

In Section \ref{sec:ktheorysoergelbimodules}, we recall basic properties of equivariant $K$-theory of flag varieties and define $K$-theory Soergel bimodules. Moreover, we give a geometric construction of morphisms of $K$-theory Soergel bimodules in terms of $K$-theoretic correspondences. This can be read independently of the other sections and does not involve any $\infty$-categories.

In Section \ref{sec:kmotivesonflagvarieties}, we discuss the category $\DK^T_{(B)}(X)$ of Bruhat-constructible $T$-equivariant $K$-motives on the flag variety and show that it can be described via chain complexes of $K$-theory Soergel bimodules.

\subsection{Acknowledgements} We thank Matthew Dyer, Marc Hoyois, Shane Kelly, Adeel Khan, Wolfgang Soergel, Catharina Stroppel, Matthias Wendt and Geordie Williamson for helpful discussions. The author was supported by Deutsche Forschungsgemeinschaft (DFG), project number 45744154, Equivariant K-motives and Koszul duality.


\section{Preliminaries on equivariant $K$-theory and $K$-motives} \label{sec:kmotives}
In this section, we define a formalism of equivariant $K$-motives $\DK^G(X)$ based on the equivariant stable motivic homotopy category introduced in \cite{hoyoisSixOperationsEquivariant2017}. Moreover, we discuss basic functorialities of $K$-theory and $G$-theory. Here, $\Lambda$ is any ring of coefficients and $k$ any base field.

\subsection{Definition}
Denote $\pt=\Spec(k).$ Let $G$ be an algebraic group over $k$ of \emph{multiplicative type}, for example, $G$ is a finite product of groups of the form $\Gm$ and $\mu_n.$ 
We use the term \emph{$G$-variety} to denote a separated $G$-scheme $X$ of finite type over $k$ which is $G$-quasi-projective, that is, admits a $G$-equivariant immersion into $\Proj(V)$ for a vector space $V$ with linear $G$-action. 
In particular, if $X$ is normal, quasi-projectivity implies $G$-quasi-projectivity.
A morphism of $G$-varieties is a morphism of schemes which is $G$-equivariant.

To any $G$-variety $X,$ \cite{hoyoisSixOperationsEquivariant2017} associates the \emph{$G$-equivariant stable motivic homotopy category} $\SH^G(X)$ which is a closed symmetric monoidal stable $\infty$-category. 
Moreover, there is a six functor formalism for $\SH^G(-)$ which fulfills properties such as base change, localisation sequences and projection formulae, see \cite[Theorem 1.1]{hoyoisSixOperationsEquivariant2017}.

In a next step, we pass from the stable homotopy category to $K$-motives. By \cite{hoyoisCdhDescentEquivariant}, for each $G$-variety $X$, there is a $E_\infty$-algebra $\KGL^G_X\in \SH^G(X)$ representing homotopy invariant $G$-equivariant $K$-theory and we define the category of \emph{$G$-equivariant $K$-motives} on $X$ as
$$\DK^G(X)\defi \operatorname{Mod}_{\KGL^G_X}(\SH^G(X))$$
the closed symmetric monoidal stable $\infty$-category of $\KGL^G_X$-modules in $\SH^G(X).$
The category of $K$-motives can be defined over any coefficient ring $\Lambda$ via
$$\DK^T(X,\Lambda)\defi \DK^T(X)\otimes_\Z \Lambda.$$
We will mostly suppress the coefficients from the notation and work with $\Lambda=\Q$ in Sections \ref{sec:stratified} and \ref{sec:kmotivesonflagvarieties}.
\subsection{Six functors}
By \cite[Theorem 1.7]{hoyoisCdhDescentEquivariant}, the collection of $E_\infty$-algebras $\KGL^G_X$ for all $G$-varieties $X$ is a \emph{cocartesian section}. That is, for a morphism $f:X\to Y$ of $G$-varieties there is a natural equivalence $f^*\KGL_Y^G\simeq \KGL_X^G$ in $\SH^G(X)$. This implies that $\DK^G(X)$ inherits the six functor formalism from $\SH^G(X)$, see \cite[Propositions 7.2.11 and 7.2.18]{cisinskiTriangulatedCategoriesMixed2019}.
We list some of the properties now, see \cite[Theorem 1.1]{hoyoisSixOperationsEquivariant2017}.
\begin{enumerate}
       \item (Pullback and pushforward) For any morphism $f:X\to Y$ of $G$-varieties there are adjoint pullback and pushforward functors
       \begin{align*}
              f^*: \DK^G(Y)&\rightleftarrows\DK^G(X): f_*.
       \end{align*}
       The functor $f^*$ is monoidal.
       \item (Exceptional pullback and pushforward) For any morphism $f:X\to Y$ of $G$-varieties there are adjoint exceptional pullback and pushforward functors
       \begin{align*}
       f_!: \DK^G(X)&\rightleftarrows\DK^G(Y): f^!.
       \end{align*}
      \item (Proper pushforward) If $f:X\to Y$ is a proper morphism of $G$-varieties there is a canonical equivalence of functors
      $$f_!\simeq f_*: \DK^G(X)\to\DK^G(Y).$$
      \item (Smooth pullback and Bott periodicity) If $f:X\to Y$ is a smooth morphism of $G$-varieties there is a canonical equivalence of functors
      $$f^!\simeq f^*:\DK^G(Y)\to\DK^G(X).$$ 
      \item (Base change) For a Cartesian square of morphism of $G$-varieties
      \[\begin{tikzcd}
             {X'} & X \\
             {Y'} & Y
             \arrow["p"', from=1-1, to=2-1]
             \arrow["q", from=1-1, to=1-2]
             \arrow["f", from=1-2, to=2-2]
             \arrow["g"', from=2-1, to=2-2]
      \end{tikzcd}\]
      there are natural equivalences of functors
      \begin{align*}
             g^*f_!&\simeq p_!q^*\text{ and }\\
             g^!f_*&\simeq p_*q^!.
      \end{align*}
      \item (Localisation) Let $j:U\hookrightarrow X \hookleftarrow X/U:i$ be a $G$-equivariant open immersion and its closed complement. Then there are homotopy cofiber sequences of functors on $\DK^G(X)$
       \begin{align*}
               j_!j^!&\to \id\to i_*i^*\text{ and }\\
              i_!i^!&\to \id\to j_*j^*.
       \end{align*}
       \item (Projection formulae) For any morphism of $G$-varieties $f$, there are natural equivalences of functors
       \begin{align*}
              f_!(-\otimes f^*(-))&\simeq f_!(-)\otimes-,\\
              \iHom(f_!(-),-)&\simeq f_*\iHom(-,f^!)\text{ and}\\
              f^!\iHom(-,-)&\simeq \iHom(f^*-,f^!-).
       \end{align*}
       \item (Homotopy invariance) If $f: E\to X$ is a G-equivariant affine bundle over a $G$-variety $X,$ then $f^*\simeq f^!:\DK^G(X)\to \DK^G(E)$ is fully faithful.
\end{enumerate}
\begin{remark}
       We note that the results of \cite{hoyoisSixOperationsEquivariant2017} work in greater generality. For example, one might work with linearly reductive groups $G.$ Moreover, the $G$-quasi-projectivity assumption can be weakened for certain nice groups $G,$ see \cite{khanGeneralizedCohomologyTheories2021}. 
\end{remark}
\begin{remark}
       A remarkable property of $K$-motives which is different from motivic sheaves or $\ell$-adic sheaves is \emph{Bott periodicity}. Namely, the reduced $K$-motive of $\Proj^1$ is isomorphic to the unit object. This implies that the Tate-twist and shift $(1)[2]$ is isomorphic to the identity in $\DK$. Bott periodicity is also reflected in the fact that $f^*\simeq f^!$ for smooth maps $f.$
\end{remark}
\subsection{$K$-motives and $K$-theory} $K$-motives compute homotopy $K$-theory and $G$-theory. In particular, by \cite[Proposition 4.6 and Remark 5.7]{hoyoisCdhDescentEquivariant} for a $G$-variety $f:X\to\pt$ we get the following equivalences of spectra
\begin{align*}
       \Map_{\DK^G(\pt)}(\Z,f_*f^*\Z)&\simeq KH^G(X)\text{ and }\\
       \Map_{\DK^G(\pt)}(\Z,f_*f^!\Z)&\simeq G^G(X).
\end{align*}
Here $KH^G(X)= \colim_{n\in\Delta^{op}} \mathbb{K}^{G}(X\times \A^n)$ denotes Weibel's \emph{homotopy $K$-theory} spectrum, see \cite[Sec. IV.12]{weibelBookIntroductionAlgebraic2013}, which is an $\A^1$-homotopy invariant version of the nonconnective $K$-theory spectrum  $\mathbb{K}$ of the category of perfect $G$-equivariant complexes on $X.$ Moreover, $G^G(X)=K(\Coh_G(X))$ denotes the $G$-equivariant $G$-theory of $X$ which is the $K$-theory spectrum of the category of $G$-equivariant coherent sheaves on $X.$

In particular, passing to homotopy groups, there are isomorphisms
\begin{align}
       \Hom_{\DK^G(\pt)}(\Z,f_*f^*\Z(p)[q])&\cong KH_{2p-q}^G(X)\text{ and }\label{eq:homindkiskh}\\
       \Hom_{\DK^G(\pt)}(\Z,f_*f^!\Z(p)[q])&\cong G_{2p-q}^G(X).
\end{align}

We note that for regular $X$ the following natural maps are equivalences of spectra
\begin{align}
       K^G(X)\to KH^G(X)\to G^G(X)\label{eq:ktheoryequalgtheory}
\end{align}
where $K^G(X)$ denotes the $K$-theory spectrum of the category of $G$-equivariant perfect complexes.

The usual functorialities of $K$-theory and $G$-theory are induced by the appropriate unit and counit maps of the adjunctions $f^*,f_*$ and $f_!,f^!$ while making use of the fact that $f_!\simeq f_*$ for $f$ proper and $f^!\simeq f^*$ for $f$ smooth. So $KH^G$ admits arbitrary pullbacks and pushfowards along smooth and proper maps, while $G^G$ admits proper pushforwards and smooth pullbacks.
\subsection{$K$-motives, Correspondences and Convolution}\label{sec:convolution} For a smooth $G$-variety $S$ and proper $G$-morphisms $p_X,p_Y,p_Z:X,Y,Z\to S$ one can use base change to identify
\begin{align}\label{eq:homsarecorrespondences}
       \Hom_{\DK^G(S)}(p_{X,!}\Z,p_{Y,!}\Z)\cong G_0^G(X\times_S Y).
\end{align}
Moreover, this identification is compatible with convolution in the following way. Consider the maps
\[\begin{tikzcd}
	{X\times_S\times Y\times Y\times_S Z} &{X\times_S\times Y\times_S Z}& {X\times_S Z.}
	\arrow["\delta"', from=1-2, to=1-1]
	\arrow["p", from=1-2, to=1-3]
\end{tikzcd}\]
Then there is a convolution product defined via
$$\star: G_0^G(X\times_SY)\times G_0^G(Y\times_SZ)\to G_0^G(X\times_S Z),\, \alpha\star\beta=p_*\delta^*(\alpha\boxtimes \beta).$$
The obvious diagram comparing composition and convolution using the isomorphisms \label{eq:homsarecorrespondences} commutes. This is shown in the context of $\DM$ and Borel--Moore motivic cohomology in \cite{fangzhouBorelMooreMotivic2016}. The same arguments apply to $\DK$ and $G$-theory.
\section{Preliminaries on stratified equivariant $K$-motives}\label{sec:stratified}
We introduce $\Ss$-stratified $G$-equivariant $K$-motives on varieties with $G$-equivariant affine stratifications and discuss basis properties, such as the existence of weight structures. In this section, we work with rational coefficients $\Lambda=\Q$ everywhere, the base field $k=\F_q$ or $k=\overline{\F}_q$ and let $\pt=\Spec(k)$.
\subsection{Constant Equivariant $K$-motives}
For an algebraic group $G$ over $k$ we denote by $R(G)=K_0(\Rep_k(G))=K_0^G(\pt)$ the representation ring. Let $G$ be an algebraic group over $k$ of multiplicative type. For a $G$-variety $X$ we consider the category of \emph{constant equivariant $K$-motives}
$$\DKT^G(X)\subset \DK^G(X)$$
as generated by the tensor unit $\Q$ by finite colimits and retracts.

In some cases, $\DKT^G(X)$ admits an explicit description in terms of modules over the representation ring $$R(G)\defi\End_{\DK^G(\pt)}(\Q)=K_0^G(\pt)=K_0(\Rep(G)).$$
\begin{proposition}\label{prop:nohigherextensionsonapoint}
    Let $G$ be a diagonalizable algebraic group and $V\in \Rep(G)$ Then
    $$\Hom_{\DK^G(V)}(\Q[n],\Q)=\left\{
        \begin{array}{ll}
        R(G) & \text{ if }n=0\text{ and } \\
        0 & \,\textrm{else.} \\
        \end{array}\right.
    $$
\end{proposition}
\begin{proof}
    By the homotopy invariance $\DK^G$ we can assume that $V=\pt$ with the trivial $G$ action. In this case, 
    $$\Hom_{\DK^G(\pt)}(\Q[n],\Q)=K_{n}^G(\pt)=K_0^G(\pt)\otimes_\Q K_{n}(\pt)$$
    where the first equality follows from \eqref{eq:homindkiskh} and \eqref{eq:ktheoryequalgtheory} and the second from \cite[Theorem 1.1(b)]{joshuaHigherKtheoryToric2015}.
    By assumption, $\pt=\Spec(\F_q)$ or $\pt=\Spec(\overline{\F}_p)$ and we use rational coefficients. Hence $K_n(\pt)=0$ for $n\neq 0$ and the statement follows.
\end{proof}
The vanishing of $\Hom_{\DKT^G(V)}(\Q[n],\Q)$ for $n>0$ allows to define the following weight structure (for an overview over weight structures and weight complex functors for $\infty$-categories, see \cite[Section 2.1.3]{eberhardtIntegralMotivicSheaves2022}) on $\DKT^G(V),$ which exists by \cite[Proposition 1.2.3(6)]{bondarkoWeightsRelativeMotives2014}.
\begin{definition}\label{def:weightstructurestratum} Let $G$ be a diagonalizable algebraic group and $V\in \Rep(G).$ The \emph{standard weight structure} $w$ on $\DKT^G(V)$ is defined as the unique weight structure on $\DKT^G(V)$ with heart $\DKT^G(V)^{w=0}$ generated by $\Q$ by finite direct sums and retracts.
\end{definition}
The vanishing of $\Hom_{\DKT^G(V)}(\Q[n],\Q)$ for $n<0$ implies that the weight complex functor 
$$\DKT^G(V)\to \Chb(\Ho(\DKT^G(V)^{w=0}))$$
to the category of chain complexes of the homotopy category of the heart of the weight structure is an equivalence of categories. The category $\Ho(\DKT^G(V)^{w=0})$ is equivalent to the category of finitely generated projective $R(G)$-modules and hence we obtain:
\begin{proposition}
    Let $G$ be a diagonalizable algebraic group and $V\in \Rep(G).$  There is an equivalence of categories between constant $G$-equivariant $K$-motives and the perfect derived category of the representation ring $R(G)$
    $$\DKT^G(V)\stackrel{\sim}{\to}\Dperf(R(G)).$$
\end{proposition}
The description is compatible with pullback/pushforward along surjective $G$-equivariant maps using the homotopy invariance of $\DK^T.$
\begin{proposition}\label{prop:surjectivelinearmap} Let $G$ be a diagonalizable algebraic group and $f:V\twoheadrightarrow W$ be an surjective map in $\Rep(G)$. Then
    $$f^*\Q\cong f^!\Q\cong \Q \in \DK^G(V) \text{ and } f_*\Q\cong f_!\Q\cong \Q \in \DK^G(W).$$
\end{proposition}
\begin{proof} Since $f$ is smooth $f^*\cong f^!$ which implies the first chain of isomorphisms. The homotopy invariance of $\DK^G$ implies the second.
\end{proof}
\begin{corollary}\label{cor:surjectivelinearmap}
    In the notation of Proposition \ref{prop:surjectivelinearmap}, the functors $f_?,f^?$ are weight exact and there are homotopy commutative diagrams
    \[\begin{tikzcd}
        {\DKT^G(V)} & {\Dperf(R(G))} \\
        {\DKT^G(W)} & {\Dperf(R(G))}
        \arrow[Rightarrow, no head, from=1-2, to=2-2]
        \arrow["{f^?}"', shift right=1, from=2-1, to=1-1]
        \arrow["{f_?}"', shift right=1, from=1-1, to=2-1]
        \arrow["\simeq", from=1-1, to=1-2]
        \arrow["\simeq", from=2-1, to=2-2]
    \end{tikzcd}\]
    for $?=*,!$ where the horizontal maps are induces from the weight complex functor.
\end{corollary}
\begin{proof}
Follows from Proposition \ref{prop:surjectivelinearmap} and the fact that the weight complex functor commutes with weight exact functors, see \cite{sosniloTheoremHeartNegative2017}.
\end{proof}
\subsection{Affine-Stratified Varieties} In this section, we consider $K$-motives for $G$-varieties with $G$-equivariant affine stratifications, that is, $G$-varieties that are stratified by $G$-represenations.
\begin{definition}\label{def:equivariantstratification}
    Let $G$ be an algebraic group and $X$ a $G$-variety. A $G$-equivariant affine stratification $\Ss$ is a decomposition
    $$X=\biguplus_{s\in \Ss}X_s$$
    of $X$ into $G$-invariant locally closed subsets, called \emph{strata}, such that for each $s\in \Ss$ the closure $\overline{X}_s$ is a union of strata and there is a $G$-equivariant isomorphism $X_s\cong V$ for some $V\in\Rep(G).$ We denote the inclusion of a stratum by $i_s:X_s\hookrightarrow X.$
\end{definition}
We need a notion of morphism between $G$-varieties with $G$-equivariant affine stratification, that is built from surjective linear maps of $G$-representations.
\begin{definition}\label{def:affinestratifiedmorphism} Let $(X,\Ss)$ and $(Y,\Ss')$ be $G$-varieties with $G$-equivariant affine stratifications. A \emph{$G$-equivariant affine stratified morphism} is a $G$-equivariant morphism $f:X\to Y$ such that 
    \begin{enumerate}
        \item for each $s\in\Ss'$, the preimage $f\inv(Y_s)$ is a union of strata.
        \item for each $X_s$ mapping into $Y_{s'}$ there is a commutative diagram
       \[\begin{tikzcd}
            {X_s} & {Y_{s'}} \\
            V & W
            \arrow["\wr", from=1-1, to=2-1]
            \arrow[two heads, from=1-1, to=1-2]
            \arrow[two heads, from=2-1, to=2-2]
            \arrow["\wr", from=1-2, to=2-2].
        \end{tikzcd}\] 
        where $V\to W$ is a surjective map in $\Rep(G).$
    \end{enumerate}
\end{definition}
We now define $K$-motives which are constant along the strata of a stratification.
\begin{definition}\label{def:constructiblekmotives}
    Let $G$ be a diagonalizable algebraic group and $(X,\Ss)$ a $G$-variety with a $G$-equivariant affine stratification. The category of \emph{$\Ss$-stratified $G$-equivariant $K$-motives on $X$} is the full subcategory
    $$\DK_\Ss^G(X)=\setbuild{M\in \DK^G(X)}{i^?_s\in \DKT^G(X_s) \text{ for } s\in \Ss, ?=*,!}.$$
\end{definition}
Next, we study well-behaved stratifications.
\begin{definition}
   In the notation of Definition \ref{def:constructiblekmotives}, the stratification is called \emph{Whitney--Tate} if
    $i_{s,?}\Q\in \DK_\Ss^G(X) \text{ for all }s\in \Ss$ and $?=+,-.$
\end{definition}
In the case of a Whitney--Tate stratification, the category $\DK_\Ss^G(X)$ is generated by the objects $i_{s,*}\Q$ (or $i_{s,!}\Q$) under finite colimits and retracts. For example, the Whitney--Tate condition is fulfilled if there are $G$-equivariant affine-stratified resolutions of stratum closures:
\begin{definition}\label{def:enoughresolutions}
        A $G$-variety $(X,\Ss)$ with a $G$-equivariant affine stratification \emph{affords  $G$-equivariant affine-stratified resolutions} if for all $s\in \Ss$ there is a $G$-equivariant map $p_s:\widetilde{X}_s\to \overline{X}_s,$ such that
    \begin{enumerate}
        \item $\widetilde{X}_s$ is smooth projective and has a $G$-equivariant affine stratification,
        \item $p_s$ is $G$-equivariant affine-stratified morphism and an isomorphism over $X_s.$
    \end{enumerate}
    
\end{definition}
There is a weight structure on constructible equivariant $K$-motives by gluing the standard weight structures on the strata, see Definition \ref{def:weightstructurestratum}.
\begin{proposition}
    Let $G$ be a diagonalizable algebraic group and $(X,\Ss)$ a $G$-variety with a Whitney--Tate $G$-equivariant affine stratification. Setting 
    \begin{align*}
        \DK_\Ss^G(X)^{w\leq 0}&=\setbuild{M\in \DK_\Ss^G(X)}{i_s^!M\in \DKT^G(X_s)^{w\leq 0}\text{ for all }s\in\Ss}\text{ and }\\
        \DK_\Ss^G(X)^{w\geq 0}&=\setbuild{M\in \DK_\Ss^G(X)}{i_s^*M\in \DKT^G(X_s)^{w\geq 0}\text{ for all }s\in\Ss}.
    \end{align*}
    defines a weight structure on $\DK_\Ss^G(X)$ that we call \emph{standard weight structure}.
\end{proposition}
\begin{proof}
    The existence follows from an iterative application of \cite[Theorem 8.2.3]{bondarkoWeightStructuresVs2010}.
\end{proof}
Stratified equivariant $K$-motives and their weight structure are compatible with  affine-stratified equivariant maps in the following way.
\begin{proposition}\label{prop:weightexactnessoffunctors}
    Let $G$ be a diagonalizable algebraic group, $(X,\Ss), (Y,\Ss')$  $G$-varieties with Whitney--Tate $G$-equivariant affine stratification and $f:X\to Y$ a $G$-equivariant affine-stratified morphism. Then the following holds.
    \begin{enumerate}
        \item The functors $f^*,f_*,f_!,f^!,\otimes$ and $\iHom$ preserve $\DK^G_\Ss.$
        \item The functors $f_!,f^*$ preserve non-negative weights.
        \item The functors $f_!,f^*$ preserve non-positive weights.
    \end{enumerate}
\end{proposition}
\begin{proof}
    Follows as in \cite[Proposition 3.8, Proposition 3.12]{eberhardtMixedMotivesGeometric2019}.
\end{proof}
The heart of the weight structure can be described in terms of the $K$-motives of resolutions of the closures of the strata.
\begin{proposition}\label{prop:weightstructuregeneratedbyresolutions}
    Let $G$ be a diagonalizable algebraic group, $(X,\Ss)$ a $G$-variety with a $G$-equivariant affine stratification which affords $G$-equivariant affine-stratified resolutions $p_s:\widetilde{X}_s\to \overline{X}_s.$ Then the heart of the weight structure $\DK_\Ss^G(X)^{w=0}$ is equal to the thick subcategory of $\DK_\Ss^G(X)$ generated by the objects $p_{s,!}\Q$ for $s\in\Ss$ by finite direct sums and retracts.
\end{proposition}
\begin{proof} By an induction on the number of strata one shows that the objects $p_{s,!}\Q$ generate the category $\DK_\Ss^G(X)$ with respect to finite colimits. By Proposition \label{prop:weightexactnessoffunctors} the objects $p_{s,!}\Q$ are contained in $\DK_\Ss^G(X)^{w=0}.$ The statement follows from the uniqueness of generated weight structures, see \cite[Proposition 1.2.3(6)]{bondarkoWeightsRelativeMotives2014}.
\end{proof}
\subsection{Pointwise Purity and Weight Complex Functor} With an additional pointwise purity assumption, stratified equivariant $K$-motives can be described in terms of their weight zero part.
\begin{definition}\label{def:pointwisepure}
    Let $G$ be a diagonalizable algebraic group, $(X,\Ss)$ be a $G$-variety with Whitney--Tate $G$-equivariant affine stratification. Let $?\in\{*,!\}.$ An object $M\in \DK_\Ss^G(X)$ is called \emph{$?$-pointwise pure} if $i_s^?M\in \DK^G(X)^{w=0}$ for all $s\in \Ss.$ The object is called \emph{pointwise pure} if it is $?$-pointwise pure for both $?=*,!.$
\end{definition}
\begin{proposition} \label{prop:homvanishingpointwisepure}
    In the notation of Definition \ref{def:pointwisepure}, let $M,N\in \DK^G(X)$ be $*$- and $!$-pointwise pure, respectively, then $\Hom_{\DK^G(X)}(M,N[n])=0$ for all $n\neq 0.$
\end{proposition}
\begin{proof}
    Follows by an induction on the number of strata and Proposition \ref{prop:nohigherextensionsonapoint}, see \cite[Lemma 3.16]{eberhardtMixedMotivesGeometric2019}.
\end{proof} 
\begin{theorem}\label{thm:pointwisepureweightcomplexfunctorequivalence}
    In the notation of Definition \ref{def:pointwisepure}, assume that all objects in $\DK^G_\Ss(X)^{w=0}$ are pointwise pure. Then the weight complex functor is an equivalence of categories
    $$\DK^G_\Ss(X)\to \Chb(\Ho\DK^G_\Ss(X)^{w=0}).$$
\end{theorem}
The assumptions of Theorem \ref{thm:pointwisepureweightcomplexfunctorequivalence} are, for example, fulfilled if there are $G$-equivariant stratified resolutions of stratum closures.
\begin{proposition}\label{prop:enoughresolutionsimplypointwisepurity}
    Under the assumptions of \ref{prop:weightstructuregeneratedbyresolutions}, all objects in $\DK^G_\Ss(X)^{w=0}$ are pointwise pure.
\end{proposition}
\begin{proof}
    The generators $p_{s,!}\Q$ of $\DK^G_\Ss(X)^{w=0}$ are pointwise pure by base change and $p_{s,!}=p_{s,*}$, see \cite[Proposition 3.15]{eberhardtMixedMotivesGeometric2019}.
\end{proof}
\section{$K$-theory Soergel bimodules}\label{sec:ktheorysoergelbimodules}
The goal of this section is to define $K$-theory Soergel bimodules. Similarly to usual, cohomological, Soergel bimodules, they arise from the equivariant $K$-theory of Bott--Samelson resolutions of Schubert varieties. We start the section with basic notations and results on representation rings and the equivariant $K$-theory of flag varieties and Bott--Samelson varieties.
Here, $\Lambda$ is any ring of coefficients and $k$ any base field.
\subsection{Flag varieties and Bott--Samelson varieties}\label{sec:flagvarandbottsamelson}
Let $G\supset B\supset T$ be a split reductive connected group over $k$ that has a simply connected derived group with a Borel subgroup $B$ and maximal torus $T$. Denote by $\Hom(T,\Gm)=X(T)\supset \Phi \supset \Phi^+$ the character lattice, set of roots and set of positive roots.
Denote by $\pW=N_G(T)/T$ the Weyl group, $\pS\subset \pW$ the set of simple reflection with respect to $B$ and by $w_0\in \pW$ the longest element. 
Let $U\subset B$ the unipotent radical, $U^-=U^{w_0}$ its opposite and $U_w=U\cap wU^{-}w^{-1}$ for $w\in W.$ 
We consider the \emph{Bruhat stratification} of the  flag variety $X=G/B$
$$X=\biguplus_{w\in \pW} X_w$$
into $B$-orbits $X_w=BwB/B$ called \emph{Bruhat cells}. For $w\in \pW$ there is a $T$-equviariant isomorphism $U_w\to X_w, u\mapsto uwB/B$ where $T$ acts on $U_w$ by conjugation and on $X_w$ by left multiplication. There is an isomorphism $U_w=\A^{\ell(w)}$ and the action of $T$ on $U_w$ is linear with set of characters $\Phi^+\cap w(\Phi^-).$ 

For a simple reflection $s\in \pS,$  let $P_s=B\cup BsB\subset G$ denote the associated parabolic subgroup. For a sequence of simple reflections $\mathbf{x}=(s_1,\dots, s_n)\in \pS^n$ denote the associated \emph{Bott-Samelson variety} and map to the flag variety by
$$p_{\mathbf{x}}:\BSvar(\mathbf{x})= P_{s_1}\times_B \dots \times_B P_{s_n}/B \to X, [p_1,\dots,p_n]\mapsto p_1\cdot \ldots \cdot p_nB/B .$$
The variety $\BSvar(\mathbf{x})$ is smooth projective and arises as quotient of $P_{s_1}\times \dots \times P_{s_n}$ by the action of $B^n$ via
$$(b_1,\dots,b_n)\cdot (p_1,\dots,p_n)=(p_1b_1\inv,b_1p_2b_2\inv,\dots,b_{n-1}p_nb_n\inv).$$
The torus $T$ acts on $\BSvar(\mathbf{x})$ from the left and there is a $T$-equivariant affine stratification on $\BSvar(\mathbf{x})$ indexed by the $2^n$ subsequences of $\mathbf{x}.$ 
Moreover, the map $p_{\mathbf{x}}$ is a $T$-equivariant affine-stratified map in the sense of Definition \ref{def:affinestratifiedmorphism}, see \cite[Proposition 2.1]{haerterichTequivariantCohomologyBottSamelson2004} and \cite[Proposition 3.0.2]{hainesProofKazhdanLusztigPurity}. 

\subsection{Representation Rings and Frobenius Extensions} 
We now discuss various representation rings and their relation to each other. The discussion mostly follows \cite{kazhdanProofDeligneLanglandsConjecture1987} and \cite{chrissRepresentationTheoryComplex2010}. We denote by $$R=R(T)=K_0^T(\pt)=\Lambda[X(T)]$$ the representation ring of $T.$ For a character $\lambda\in X(T)$, we write $e^\lambda$ for the corresponding element in $R.$ This way, $e^\lambda=[k_\lambda]$ denotes the class of the one-dimensional representation $k_\lambda$ on which $T$ acts via $\lambda$.
The ring $R(T)$ is isomorphic to a Laurent polynomial ring in $\operatorname{rank}(T)$ many variables. Moreover, there is a natural action of $\pW$ on $R(T).$

The representation ring $R(G)=K^G_0(\pt)$ is related to $R=R(T)$ via two natural maps
$$\Ind_T^G: R(T)\leftrightarrows R(G):\Res_G^T.$$
We will describe these maps in detail and see that they form a \emph{Frobenius extension}. 

The map $\Res_G^T$ is an injective algebra homomorphism defined by restricting a $G$-representation to $T.$ 
The image of $\Res_G^T$ are exactly the $\pW$-invariants and we hence identify $R(G)=R^{\pW}\subset R.$

The map $\Ind_T^G$ is obtained by inducing representations from $T$ to $G$. 
It maps the class $[V]\in R(T)$ of a representation $V$ of $T$ to the alternating sum of the cohomology groups of the $G$-equivariant vector bundle $G\times_BV$ on the flag variety $X$, where $B$ acts on $V$ via the quotient map $B\to T.$ Namely, $\Ind_T^G$ is the composition of maps
\[\begin{tikzcd}[column sep=small, row sep=small]
	{R=R(T)} & {K_0^G(X)} & {K^G_0(\pt)=R(G)=R^\pW} \\
	{[V]} & {[G\times_BV]} & {\sum_i(-1)^i[H^i(X,G\times_BV)]}
	\arrow["\sim", from=1-1, to=1-2]
	\arrow["{\pi_*}", from=1-2, to=1-3]
	\arrow[maps to, from=2-1, to=2-2]
	\arrow[maps to, from=2-2, to=2-3]
\end{tikzcd}\]
where the first isomorphism comes from the induction equivalence
$$R(T)=K_0^T(\pt)\cong K_0^B(\pt)\cong K_0^G(G\times_B\pt)= K_0^G(X)$$
and $\pi_*$ is pushforward along the projection map $\pi: X\to \pt.$ 

The map $\Ind_T^G$ is an $R(G)$-module homomorphism. Namely, let $[V']\in R(G)$ be the class of a representation of $G$. Then there is a $G$-equivariant trivialisation
$$G\times_B V' \stackrel{\sim}{\to} G/B\times V', [g,v]\mapsto (gB, gv)$$
which shows that $\Ind_T^G(\Res_G^T[V'])=[V'].$ Similarly, if $[V]\in R(T)$ is the class of a representation $V$ of $T,$ then $G\times_B(V'\otimes V)$ is the tensor product of the vector bundles $G\times_BV'\cong G/B\times V'$ and $G\times_BV.$ It follows that  $\Ind_T^G(\Res_G^T([V'])[V])=[V']\Ind_T^G([V]).$

The \emph{Weyl character formula} allows to explicitly compute $\Ind_T^G$ as 
$$\Ind_T^G(e^\lambda)=\frac{\sum_{w\in \pW}(-1)^{\ell(w)}e^{w\cdot\lambda}}{\prod_{\alpha\in \Phi^+}(1-e^{-\alpha})}$$
where $w\cdot\lambda=w(\lambda+\rho)-\rho$ denotes the dot-action of $\pW$ and $\rho=\frac{1}{2}\sum_{\alpha\in\Phi^+}\alpha$ is the half-sum of all positive roots.

The map $\Ind_T^G$ induces a pairing
$$\langle\,,\,\rangle: R(T)\otimes R(T)\to R(G), \langle[V],[V']\rangle=\Ind_T^G([V][V'])=\pi_*([G\times_B(V\otimes V')]).$$
There is an $R(G)$-basis $\{e_w\}_{w\in \pW}$ of $R(T)$ constructed in \cite{steinbergTheoremPittie1975} such that 
$$\det(\langle e_w,e_{w'} \rangle)_{w,w'}=1,$$
see \cite[Proposition 1.6]{kazhdanProofDeligneLanglandsConjecture1987}. 
This implies that there is a dual basis $\{e^*_w\}_{w\in \pW}$ such that $\langle e_w,e^*_{w'} \rangle=\delta_{w,w'}.$ 
 Hence $\Res_G^T$ and $\Ind_T^G$ form a Frobenius extension. It follows that the functors
$$R(T)\otimes_{R(G)}-: \Mod_{R(G)}\leftrightarrows \Mod_{R(T)}:\Hom_{R(T)}(R(T),-)$$
are adjoint in both ways.

We remark that the discussion also applies to standard parabolic subgroups $B\subset P\subset G$ by taking $G=L=P/\operatorname{Rad}_u(P)$ as Levi factor of $P$ and using that $R(P)\cong R(L)$ and $P/B\cong L/(B\cap L).$
\subsection{The rank two case} The previous discussion specializes to the following formulas for minimal parabolics $P_s=B\cup BsB\subset G$ for simple reflections $s\in S$. Namely, one can identify $R(P_s)=R^s$ and there is a Frobenius extension
$$\Ind_T^{P_s}: R(T)\leftrightarrows R(P):\Res_{P_s}^T$$
where $\Res_{P_s}^T$ corresponds to the inclusion $R^s\subset R$ and $\Ind_T^{P_s}$ is given by
$$\Ind_T^{P_s}(e^\lambda)=\frac{e^\lambda-e^{s\cdot\lambda}}{1-e^{-\alpha_s}}=\frac{e^{\lambda+\alpha_s/2}-e^{s(\lambda)-\alpha_s/2}}{e^{\alpha_s/2}-e^{-\alpha_s/2}}$$
where $\alpha_s\in \Phi^+$ is the simple root corresponding to $s.$ 
Hence, the $\Delta_s=\Ind_T^{P_s}\Res_{P_s}^T$ for $s\in S$ are the \emph{Demazure operators} on $R=R(T)$, see \cite{demazureNouvelleFormuleCaracteres1975}.

\subsection{Equivariant $K$-theory of Flag and Bott--Samelson varieties}
We study the $T$-equivariant $K$-theory of the flag variety, Bruhat cells and Bott-Samelson varieties. 
There are isomorphisms
\begin{align}
	K_0^T(G/B)\cong K^{T\times T}_0(G)\cong R\otimes_{R^{\pW}}R \label{eq:ktheoryflagvariety}
\end{align}
where the second isomorphism is induced from the pullback $R\otimes R=K^{T\times T}(\pt)\to K^{T\times T}(G)$. In particular, we can interpret modules over $K_0^T(G/B)\cong R\otimes_{R^{\pW}}R$ as $R$-bimodules.

For a stratum $X_w=BwB/B\subset X$ one has
\begin{align}
	K_0^T(X_w)\cong K_0^{T\times T}(wT)\cong R_w \label{eq:ktheorybruhatcell}
\end{align}
where $R_w$ is isomorphic to $R$ as a ring but has a twisted $R$-bimodule structure, given by $r\cdot m=rm$ and $m\cdot r =mw(r)$ for $r\in R, m\in R_w.$

Next, we compute the $T$-equivariant $K$-theory of Bott--Samelson varieties. For this, we make use of the following statement:
\begin{lemma}\label{lem:balancedproductwithparabolic}
	Let $B\subset P\subset G$ be a standard parabolic. Let $X$ be a $B$-variety. 
	Then there is a natural isomorphism
	$$K_0^T(P\times_BY)\cong R\otimes_{R^{W_P}}K_0^T(Y).$$
\end{lemma}
\begin{proof} There is the following chain of isomorphisms
	\begin{align*}
		K_0^T(P\times_BX)&\cong R\otimes_{R^{W_P}}K_0^P(P\times_B X)\\
		&\cong  R \otimes_{R^{W_P}}K_0^B(X)\\
		&\cong R\otimes_{R^{W_P}}K_0^T(X).
	\end{align*}
	The first isomorphism is \cite[Theorem 6.1.22]{chrissRepresentationTheoryComplex2010}, the second the induction equivalence and the last the reduction property.
\end{proof}
Let $\mathbf{x}=(s_1,\dots, s_n)\in \pS^n$ be a sequence of simple reflections. By applying Lemma \ref{lem:balancedproductwithparabolic} inductively, one obtains that
	\begin{align}
		K_0^T(\BSvar(\mathbf{x}))=K_0^T(P_{s_1}\times_B \dots \times_B P_{s_n}/B)\cong R\otimes_{R^{s_1}} \otimes \dots \otimes_{R^{s_n}} R\label{eq:ktheorybottsamelson}
	\end{align}
as an $R$-bimodule.
\subsection{$K$-theory Soergel bimodules}
Soergel bimodules arise from (direct summands of) the $T$-equivariant cohomology of Bott--Samelson varieties, interpreted as bimodules over the $T$-equivariant cohomology ring of a point $H^\bullet_T(\pt)=H^\bullet(BT).$
It is hence natural to define $K$-theory Soergel bimodules in the same way, replacing equivariant cohomology by equivariant $K$-theory.
\begin{definition}
	The category of \emph{$K$-theory Soergel bimodules} $\SBim_R$ is the full thick subcategory of the category of $\Ho\Bim_R$ generated by the $R$-bimodules 
	$$K_0^T(\BSvar(\mathbf{x}))=R\otimes_{R^{s_1}} \otimes \dots \otimes_{R^{s_n}}R$$
	for all sequences $\mathbf{x}=(s_1,\dots, s_n)\in \pS^n$ of simple reflections by finite direct sums and retracts.
\end{definition}
\begin{remark}
		In fact, it will turn out that (with rational coefficients) the category is $\SBim_R$ is already generated by the collection of bimodules $K_0^T(\BSvar(\underline{w}))$ for any fixed choice of reduced expressions $\underline{w}$ for the elements $w\in \pW.$ This follows from the geometric description in terms of weight zero $K$-motives, see Corollary \ref{cor:purekmotivesviasoergelbimodules}, and Proposition \ref{prop:weightexactnessoffunctors}.
\end{remark}
\subsection{$K$-theory Soergel bimodules via Convolution} We will now show how homomorphisms between $K$-theory Soergel bimodules can be described via a convolution product. This yields an equivalent definition of the category $K$-theory Soergel bimodules via correspondences.

Namely, for two sequences of simple reflections $\mathbf{x}\in \pS^n$ and $\mathbf{y}\in \pS^m$ convolution defines a natural map
$$\act:G_0^T(\BSvar(\mathbf{x})\times_X \BSvar(\mathbf{y}))\to \Hom_{K_0^T(X)}(K_0^T(\BSvar(\mathbf{x})), K_0^T(\BSvar(\mathbf{y}))), \beta\mapsto (\alpha\mapsto \alpha\star \beta)$$
where $\alpha\star \beta=p_*\delta^*(\alpha\boxtimes \beta)$ is the convolution of $\alpha$ and $\beta$ and the maps $\delta^*$ and $p_*$  are induced by the diagonal and projection maps
\[\begin{tikzcd}
	{\BSvar(\mathbf{x})\times\BSvar(\mathbf{x})\times_X\BSvar(\mathbf{y}) } & {\BSvar(\mathbf{x})\times_X \BSvar(\mathbf{y})} & {\BSvar(\mathbf{y}),}
	\arrow["\delta"', from=1-2, to=1-1]
	\arrow["p", from=1-2, to=1-3]
\end{tikzcd}\]
see Section \ref{sec:convolution}.
\begin{theorem}\label{thm:erweiterungssatzconvolutionktheory} The map $\act$ is an isomorphism.
\end{theorem}
\begin{proof}\emph{Step 1:}
	We reduce the statement to the case when $\mathbf{x}$ is the empty sequence and hence $\BSvar(\mathbf{x})=B/B=X_e.$

	For this, let $s\in \pS$ be a simple reflection and write $s\mathbf{x}\in \pS^{n+1}$ for the concatenation. Then $\BSvar(\mathbf{sx})=P_s\times_B\BSvar(\mathbf{x}).$ We abbreviate $P=P_s,$ $M=\BSvar(\mathbf{x})$ and $N=\BSvar(\mathbf{y}).$ Our goal is to construct a commutative diagram
\[		\begin{tikzcd}
		{G_0^B((P\times_BM)\times_X N)} & {G_0^B(M\times_X (P\times_B N))} \\
		{\Hom_{R\otimes R}(K_0^B(P\times_B M), K_0^B(N))} & {\Hom_{R\otimes R}(K_0^B(M),K_0^B(P\times_BN))}
		\arrow["\sim", from=1-1, to=1-2]
		\arrow["\act"', from=1-1, to=2-1]
		\arrow["\act", from=1-2, to=2-2]
		\arrow["\sim", from=2-1, to=2-2]
	\end{tikzcd}\]
	The upper horizontal isomorphism can be constructed as follows. The isomorphism
	$$\phi: (P\times M)\times_X N\to  M\times_X (P\times N), (p,m,n)\mapsto (m,p^{-1},n)$$
	is equivariant with respect to the actions by $B\times B$ on $(P\times M)\times_X$ and $M\times_X (P\times N)$ given by 
	$(b_1,b_2)(p,m,n)=(b_1pb_2^{-1},b_2m,b_1n)$ and $(b_1,b_2)(m,p,n)=(b_2m,b_2pb_1^{-1},b_1n),$ respectively. 
	Hence, $\phi$ induces via the induction equivalence an isomorphism
	\begin{align*}
	G_0^B((P\times_BM)\times_X N)&= G_0^{B\times B}((P\times M)\times_X N) \\\stackrel{\sim}{\to}G_0^{B\times B}(M\times_X (P\times N))&=G_0^B(M\times_X (P\times_B N)).
	\end{align*}
	The lower horizontal isomorphism in the commutative diagram comes from the following chain of isomorphisms
	\begin{align*}
		\Hom_{R\otimes R}(K_0^B(P\times_B M), K_0^B(N))&\stackrel{\sim}{\to}\Hom_{R\otimes R}(R\otimes_{R^s}K_0^B(M), K_0^B(N))\\
		&\stackrel{\sim}{\to} \Hom_{R^s\otimes R}(K_0^B(M), K_0^B(N))\\
		&\stackrel{\sim}{\leftarrow}\Hom_{R\otimes R}(K_0^B(M), R\otimes_{R^s}K_0^B(N))\\
		&\stackrel{\sim}{\leftarrow}\Hom_{R\otimes R}(K_0^B(M),K_0^B(P\times_BN)).
	\end{align*}
	The first and last isomorphisms are given by Lemma \ref{lem:balancedproductwithparabolic}. The second isomorphism is the Hom-tensor adjunction and sends a map $f$ to $(x\mapsto f(1\otimes x)).$ The third isomorphism comes from the
	Frobenius extension $R^s\subset R$ with trace map $\Delta_s:R\to R_s$ and is given by pushforward along the map $r\otimes y\mapsto \Delta_s(r)y.$ In total, the composition of the first two isomorphisms is induced by the pullback along the map $P\times M\to M.$ Dually, the composition of the last two isomorphism is induced by the pushforward along the map $P\times N\to N.$

	The commutativity of the above square boils down to the commutativity of the diagram
\[\begin{tikzcd}[column sep=-2pt]
	{G_0^{B\times B}((P\times M)\times_X N)} && {G_0^{B\times B}(M\times_X (P\times N))} \\
	& {\Hom(K_0^B(M), K_0^B(N))}
	\arrow["{\act'}"', from=1-1, to=2-2]
	\arrow["{\act''}", from=1-3, to=2-2]
	\arrow["{\phi^*}", from=1-1, to=1-3]
\end{tikzcd}\]
where $\act'$ and $\act''$ are defined via the exterior tensor product as well as pullback and pushforward along the diagrams
\[\begin{tikzcd}[column sep=15pt]
	{M\times (P\times M)\times_XN} & {(P\times M)\times (P\times M)\times_XN} & {(P\times M)\times_XN} & N
	\arrow["\pi"', from=1-2, to=1-1]
	\arrow["\delta"', from=1-3, to=1-2]
	\arrow["p", from=1-3, to=1-4]
\end{tikzcd}\]
and
\[\begin{tikzcd}
	{M\times M\times_X(P\times N)} & {  M\times_X(P\times N)} & {P\times N} & N,
	\arrow["\delta"', from=1-2, to=1-1]
	\arrow["\pi", from=1-3, to=1-4]
	\arrow["p", from=1-2, to=1-3]
\end{tikzcd}
\]
respectively. We remark that the pushforwards are well-defined, since they can be represented by pushforwards along proper maps when taking the quotient by the appropriate free $B$-action. The two actions clearly agree with respect to the map $\phi.$

\emph{Step 2:} By the first step, it suffices to show that the map
$$\act:G_0^T(X_e \times_X \BSvar(\mathbf{y}))\to \Hom_{K_0^T(X)}(K_0^T(X_e), K_0^T(\BSvar(\mathbf{y})))$$
is an isomorphism. To see this, we abbreviate $N=\BSvar(\mathbf{y})$ and $N_w=X_w\times_X N.$ 
Denote by $i: N_e\to N\leftarrow N\backslash N_e:u$ the inclusions. We identify $K_0^T(X_e)=R$. Then $\act(\alpha)(1)=i_*\alpha.$
Each space $N_w$ admits a stratification such that the strata are affine bundles over $X_w.$ By the cellular fibration lemma, see \cite[Lemma 5.5.1]{chrissRepresentationTheoryComplex2010}, it follows that $K_0^T(N_w)=G_0^T(N_w)$ admits a filtration with subquotients of the form $K_0^T(X_w)=R_w.$ Moreover, there is a short exact sequence \[\begin{tikzcd}[column sep=10pt]
	0 & {K_0^T(N_e)=G_0^T(N_e)} & {G_0^T(N)} & {G_0^T(N\backslash N_e)} & 0
	\arrow[from=1-1, to=1-2]
	\arrow["{i_*}", from=1-2, to=1-3]
	\arrow["{u^*}", from=1-3, to=1-4]
	\arrow[from=1-4, to=1-5]
\end{tikzcd}\]
where $G_0^T(N\backslash N_e)$ is a succesive extension of modules of the form $R_w$ for $w\neq e.$
In the associated exact sequence
\[\begin{tikzcd}[column sep=10pt]
	0 & {\Hom(R,G_0^T(N_e))} & {\Hom(R,G_0^T(N))} & {\Hom(R,G_0^T(N\backslash N_e))}
	\arrow[from=1-1, to=1-2]
	\arrow[from=1-2, to=1-3]
	\arrow[from=1-3, to=1-4]
\end{tikzcd}\]
where we abbreviate $\Hom=\Hom_{K_0^T(X)}$ the right hand term vanishes since $$\Hom_{K_0^T(X)}(R,R_w)=\Hom_{R\otimes R}(R,R_w)=0$$ for $w\neq e.$ This implies that $\act$ is an isomorphism.
\end{proof}
\begin{remark}\label{rem:remafterweiterungssatz}
\begin{enumerate}[wide, labelwidth=!, labelindent=0pt]
	\item The isomorphism $\act$ is compatible with composition in the following sense. If $\mathbf z\in \pS^k$ is a third sequence of simple reflections, one can define the convolution product
	$$\star:G_0^T(\BSvar(\mathbf{x})\times_X \BSvar(\mathbf{y}))\times G_0^T(\BSvar(\mathbf{y})\times_X \BSvar(\mathbf{z}))\to G_0^T(\BSvar(\mathbf{x})\times_X \BSvar(\mathbf{z})),$$
	see Section \ref{sec:convolution}.
	By associativity of convolution $\act(\beta)\circ \act(\alpha) = \act(\alpha\star\beta).$
	\item The above discussion yields the following equivalent construction of the category $\SBim_R.$ Namely, consider the category of $K$-theoretic correspondences of Bott--Samelson resolutions $\KCorr$ with objects sequences of simple reflections $\mathbf{x}\in S^n$ and morphisms given by $\Hom_{\KCorr}(\mathbf{x},\mathbf{y})=K_0^T(\BSvar(\mathbf{y})\times_X \BSvar(\mathbf{x}))$ and composition given by convolution $\star.$ Then the maps $\act$ define a functor $\KCorr\to \SBim_R$ which is fully faithful by Theorem \ref{thm:erweiterungssatzconvolutionktheory}. In fact, the induced functor from the Karoubian envelope $\Kar(\KCorr)\to \SBim_R$ yields an equivalence of categories.
	\item The category $\Kar(\KCorr)$ has a more conceptual construction. Namely, there is an equivalence $\DK_{(B)}^T(X)^{w=0}\stackrel{\sim}{\to} \Kar(\KCorr)$ with the category of weight zero objects in the category of Bruhat-stratified $T$-equivariant $K$-motives on the flag variety. In this context, Theorem \ref{thm:erweiterungssatzconvolutionktheory} can be seen as a $K$-theoretic analogue of Soergel's Erweiterungssatz. It is equivalent to the statement that the functor
	$\Kyp: \DK_{(B)}^T(X)^{w=0}\to \SBim_R,$ which sends a $K$-motive to its $K$-theory, see Remark \ref{rem:remarkaftererweiterungssatzforkmotives}, is fully faithful. In fact, our proof closely follows the proof of the Erweiterungssatz in the context of equivariant motives, see \cite[Proposition III.6.11]{soergelEquivariantMotivesGeometric2018}.
\end{enumerate}
\end{remark}

\section{$K$-theory Soergel bimodules via $K$-motives on flag varieties}\label{sec:kmotivesonflagvarieties}
We now combine the results from Section \ref{sec:stratified} and \ref{sec:ktheorysoergelbimodules} to obtain a combinatorial description of Bruhat-stratified torus-equivariant $K$-motives on flag varieties in terms of (complexes of) $K$-theory Soergel bimodules. In this section, our ring of coefficients is $\Q$ and $k=\F_q$ or $\overline{\F}_p.$
\subsection{Bruhat-stratified $K$-motives} We continue in the notation of \ref{sec:flagvarandbottsamelson}. We consider the flag variety $X=G/B$ with its action by the maximal torus $T.$ By the discussion there, the Bruhat stratification is a $T$-equivariant affine stratification of $X$ in the sense of Definition \ref{def:equivariantstratification} and we denote it by $(B).$ It hence makes sense to consider the category $\DK_{(B)}^T(X)$ of \emph{Bruhat-stratified $T$-equivariant $K$-motives on the flag variety}. 

Moreover, for a reduced expression $\underline{w}=(s_1,\dots,s_n)$ of an element $w\in \pW,$ the map  $p_{\underline{w}}:\BSvar(\underline{w})\to X$ provides a resolution of singularities of the \emph{Schubert variety} $\overline{X}_w$ and hence $X$ affords $T$-equivariant affine-stratified resolutions in the sense of Definition \ref{def:enoughresolutions}. This shows that the Bruhat-stratification is Whitney--Tate and that there is a weight stucture on  $\DK_{(B)}^T(X)$ such that the objects in the heart $\DK^T_{(B)}(X)^{w=0})$ are pointwise pure by Proposition \ref{prop:enoughresolutionsimplypointwisepurity}. Hence, Theorem \ref{thm:pointwisepureweightcomplexfunctorequivalence} implies the following:
\begin{theorem}\label{thm:weightcomplexfunctorflagvariety}
	The weight complex functor induces an equivalence of categories
	$$\DK^T_{(B)}(X)\to \Chb(\Ho\DK^T_{(B)}(X)^{w=0}).$$
\end{theorem}
\subsection{A combinatorial description}\label{sec:combinatorialdescriptionofkmotives}

Let $\mathbf{x}\in \pS^n,\mathbf{y}\in \pS^m$ be sequences of simple reflections and denote by $p_\mathbf{x}:\BSvar(\mathbf{x})\to X$ and $p_\mathbf{y}:\BSvar(\mathbf{y})\to X$ the Bott--Samelson resolutions. Then combining the discussion in Section \ref{sec:convolution} and Theorem \ref{thm:erweiterungssatzconvolutionktheory}
we obtain isomorphisms
\begin{align*}
	\Hom_{\DK^T(X)}(p_{\mathbf{x},!}\Q,p_{\mathbf{y},!}\Q)&\stackrel{\sim}{\to} G^T_0(\BSvar(\mathbf{x})\times_X\BSvar(\mathbf{y}))\\
	&\stackrel{\sim}{\to} \Hom_{K_0^T(X)}(K_0^T(\BSvar(\mathbf{x})),K_0^T(\BSvar(\mathbf{y}))).
\end{align*}
compatible with composition.
Since the categories $\Ho\DK^T_{(B)}(X)^{w=0}$ and $\SBim_R$ are generated by direct sums and direct summands of the objects $p_{\mathbf{x},!}\Q$ and $K_0^T(\BSvar(\mathbf{x})),$ respectively, we obtain the following statement.
\begin{corollary}\label{cor:purekmotivesviasoergelbimodules}
	There is an equivalence of categories
	\begin{align*}
		\Ho\DK^T_{(B)}(X)^{w=0}&\stackrel{\sim}{\to}\SBim_R.
	\end{align*}
\end{corollary}
Together with Theorem \ref{thm:weightcomplexfunctorflagvariety} this yields:
\begin{corollary}\label{cor:kmotivesviacomplexesofsoergelbimodules}
	There is an equivalence of categories
	\begin{align*}
		\DK^T_{(B)}(X)&\stackrel{\sim}{\to}\Chb(\SBim_R).
	\end{align*}
\end{corollary}
 \begin{remark}\label{rem:remarkaftererweiterungssatzforkmotives}
		The equivalence $\Ho\DK^T_{(B)}(X)^{w=0}\stackrel{\sim}{\to}\SBim_R$ can also be constructed via the functor
		$$\Kyp: \Ho\DK^T(X)\to \Mod_{K^T(X)},\, M\mapsto \Hom_{\DK^T(X)}(\Q, M).$$
		Hence, Corollary \label{cor:purekmotivesviasoergelbimodules} can be seen as a $K$-theoretic analog of Soergel's Erweiterungssatz, see also Remark \ref{rem:remafterweiterungssatz}.
 \end{remark}

\bibliographystyle{amsalpha} 
\bibliography{main}

\end{document}